# Adaptive estimation of linear functionals by model selection

**Béatrice Laurent**

*Institut de Mathématiques (UMR 5219), INSA de Toulouse, Université de Toulouse, France*
*e-mail:* beatrice.laurent@insa-toulouse.fr

**Carenne Ludeña**

*IVIC, Venezuela*
*e-mail:* cludena@euclides.ivic.ve

**Clémentine Prieur**

*Institut de Mathématiques (UMR 5219), INSA de Toulouse, Université de Toulouse, France*
*e-mail:* clementine.prieur@insa-toulouse.fr

**Abstract:** We propose an estimation procedure for linear functionals based on Gaussian model selection techniques. We show that the procedure is adaptive, and we give a non asymptotic oracle inequality for the risk of the selected estimator with respect to the $\mathbb{L}_p$ loss. An application to the problem of estimating a signal or its $r^{th}$ derivative at a given point is developed and minimax rates are proved to hold uniformly over Besov balls. We also apply our non asymptotic oracle inequality to the estimation of the mean of the signal on an interval with length depending on the noise level. Simulations are included to illustrate the performances of the procedure for the estimation of a function at a given point. Our method provides a pointwise adaptive estimator.

## 1. Introduction

We consider the following model:

$$Y(t) = \langle s, t \rangle + \frac{\sigma}{\sqrt{n}} L(t), \text{ for all } t \in \mathbb{H}, \quad (1.1)$$

where $\mathbb{H}$ is a separable Hilbert space endowed with the scalar product $\langle ., . \rangle$ and $L$ is some centered Gaussian isonormal process, which means that $L$ maps isometrically $\mathbb{H}$ onto some Gaussian subspace of $\mathbb{L}_2(\Omega)$, where $(\Omega, \mathcal{G}, P)$ is some canonical probability space. This framework includes the finite dimensional Gaussian regression model, the Gaussian sequence model and the multivariate white noise model.

Let $T$ be a linear functional defined over $\mathcal{S} \subset \mathbb{H}$. In this paper we consider the problem of estimating $T(s)$, based on the observation of $(Y(t), t \in \mathbb{H})$. Our main goal will be to develop procedures which adapt to the smoothness of the underlying function $s$ in the framework of model selection as proposed by Barron et al. [3].

Minimax theory for estimating linear functionals is well developed in the Gaussian setting. Ibragimov and Hasminskii [16] obtained the best minimax linear estimator over classes of smooth functions. For any convex parameter space $\mathcal{F}$, the minimax mean squared error is of the order of the modulus of continuity of the functional over $\mathcal{F}$ (see Donoho and Liu [13], Donoho [12] and Cai and Low [8, 9], the latter being a generalization to certain non convex parameter spaces). These authors have also constructed procedures which have maximum risk close, up to a small factor, to the minimax rates. However, these rates cannot be attained when dealing with adaptive estimation over several classes of parameters. For the problem of estimating a function at a given point, Lepski [19] showed that it is necessary to include a logarithmic factor in the mean squared error when dealing simultaneously with two Lipschitz classes. For general parameter spaces, Cai and Low [11, 10] show that it is necessary to include a between class modulus of continuity to quantify precisely the degree of adaptability for the estimation of a linear functional with respect to the mean squared error. They also proposed an adaptive estimator based on multiple tests over an ordered sequence of parameter spaces. Their methodology thus resembles "Lepski's method" (see for example [19, 20, 21]) in the sense that the estimation procedure chooses the best possible over a finite selection of parameter spaces. This point of view is also developed in Klemelä and Tsybakov [17], where the authors construct an asymptotically sharp adaptive estimator of $T(s)$ based on kernel methods. They assume that the signal $s$ belongs to a class of regular functions, the index of regularity being bounded from above and below by known constants. Lepski and Spokoiny [23], Lepski, Mammen and Spokoiny [22] propose methods based on kernel estimates with variable bandwidth selector for pointwise adaptive estimation in the Gaussian white noise model.

Model selection methods for adaptive estimation have been initiated in a series of papers by Birgé and Massart (see for example Birgé and Massart [5, 6], Barron, Birgé and Massart [3]). These methods have been used in the framework of the regression with fixed or random design, to estimate the regression function by Baraud [1] and [2]. In this article, following Birgé [4] we take a model selection point of view at adaptive estimation via Lepski's method. In order to construct the adaptive estimator of the linear functional we shall choose among an ordered family of finite dimensional linear subspaces of $\mathbb{H}$. Over each subspace we consider an estimator based on projection methods and the problem is thus establishing a best possible procedure for determining the subspace. The main issue here is that, unlike the case of penalized least squares, the bias of the estimator is not a monotonically decreasing sequence over the family of nested subspaces. Hence it is necessary to modify the procedure as developed by Birgé [4] in order to obtain an appropriate estimator of the bias. The main advantage of our formulation is that it allows to obtain non asymptotic oracle inequalities for general linear functionals. In the framework of the Gaussian sequence model

$$Y_i = \theta_i + \sigma \xi_i, i \in \mathbb{N}$$

where the $\xi'_i s$ are i.i.d. standard Gaussian variables, Golubev and Levit [14] obtain an oracle inequality for the estimation of a general linear functional of $\theta = (\theta_i)_{i \in \mathbb{N}}$. They assume that the $\theta'_i s$ are independent centered Gaussian variables, this is not the framework that we consider in the present paper.

We propose a general method to estimate a linear functional of $s$ in the framework of Model (1.1). We apply this general procedure to estimate the value of the $r^{th}$ derivative of a function at a point. For this problem, we provide minimax rates uniformly over Besov balls that correspond to the rates established by Lepski [19]. We also give an application of our procedure to pointwise adaptive estimation in a multidimensional framework. Moreover, since we have obtained a non asymptotic oracle inequality, we are able to apply our result to the estimation of linear functionals that depend on the noise level (or on the number of observations). In the white noise model, we consider the estimation of the mean of the signal on an interval with length depending on the noise level. The interesting fact in this case is that we obtain two kinds of rates of convergence, according to the relationship between the length of the interval, the noise level, and the regularity of the signal. When the length of the interval is too small, this problem is as hard as estimating the signal at some fixed point and when the length of the interval is large, the functional can be estimated at the parametric rate $1/\sqrt{n}$. All intermediate rates are obtained as the length of the interval grows. We present simulation results to estimate a function at a point, and we compare our method to a global (not pointwise) model selection estimator and to an estimator based on wavelet shrinkage. Our method provides a locally adaptive estimator of a regression function $s$ on $[0,1]$. It has good properties when estimating functions that are very oscillating over some regions and nearly flat over other ones.

The article is organized as follows. In Section 2 we present the framework, the estimation procedure and our main result. In Section 3 we develop three examples: estimating the value of the $r^{th}$ derivative of a function at a point using a multiresolution analysis, estimating the mean of the signal on an interval with length depending on the noise level and estimating the value of a multidimensional function at a point. In Section 4 we present the simulation study. Proofs of our main results are given in Section 5.

## 2. Main results

### *2.1. The framework*

Given some separable Hilbert space $\mathbb{H}$, endowed with the scalar product $\langle .,. \rangle$, one observes $(Y(t), t \in \mathbb{H})$ as defined by Model (1.1). Since $L$ is some centered Gaussian isonormal process defined on $\mathbb{H}$, we have that for all $t, u \in \mathbb{H}$, $\text{Cov}(L(t), L(u)) = \langle t, u \rangle$.

Let us first consider three particular cases of Model (1.1).

### The finite dimensional Gaussian regression.
One observes
$$Y_i = s_i + \varepsilon_i, i = 1 \ldots, n$$
where $(\varepsilon_1, \ldots, \varepsilon_n)$ are independent standard normal variables.

We consider $\mathbb{H} = \mathbb{R}^n$ endowed with the scalar product $\langle x, y \rangle = \frac{1}{n} \sum_{i=1}^n x_i y_i$ and set $s = (s_1, \ldots, s_n)$.

Model (1.1) is obtained by setting, for all $t = (t_1, \ldots, t_n) \in \mathbb{R}^n$, $Y(t) = \frac{1}{n} \sum_{i=1}^n t_i Y_i$ and $L(t) = \frac{1}{\sqrt{n}} \sum_{i=1}^n t_i \varepsilon_i$.

### The Gaussian sequence model.
In the Gaussian sequence model, one observes
$$Y_\lambda = \beta_\lambda + \frac{1}{\sqrt{n}} \varepsilon_\lambda, \lambda \in \mathbb{N}^*, \tag{2.2}$$
where $(\varepsilon_\lambda)_{\lambda \in \mathbb{N}^*}$ is a sequence of independent standard normal variables.

Setting $\mathbb{H} = l_2(\mathbb{N}^*)$ endowed with the usual scalar product $\langle \beta, \gamma \rangle = \sum_{\lambda \in \mathbb{N}^*} \beta_\lambda \gamma_\lambda$ and $s = (\beta_\lambda)_{\lambda \in \mathbb{N}^*}$, we define for any $t = (\alpha_\lambda)_{\lambda \in \mathbb{N}^*} \in \mathbb{H}$, $Y(t) = \sum_{\lambda \in \mathbb{N}^*} \alpha_\lambda Y_\lambda$ and $L(t) = \sum_{\lambda \in \mathbb{N}^*} \alpha_\lambda \varepsilon_\lambda$ and we see that (2.2) is a particular case of Model (1.1).

### The multivariate white noise model.
One observes
$$Z(x) = \int_{[0,1]^d} \mathbb{1}_{[0,x_1] \times \cdots \times [0,x_d]}(u) s(u) du + \frac{1}{\sqrt{n}} W(x)$$
for all $x = (x_1, \ldots, x_d) \in [0,1]^d$, where $W$ is the standard Wiener Process on $[0,1]^d$. We consider $\mathbb{H} = \mathbb{L}_2([0,1]^d)$ endowed with its usual scalar product.

We set $Y(t) = \int_{[0,1]^d} t(u) dZ(u)$ and $L(t) = \int_{[0,1]^d} t(u) dW(u)$.

Our purpose is to propose new adaptive estimators of $T(s)$, where $T$ is a linear functional, from observation (1.1).

### *2.2. The estimation procedure*

We consider a finite or countable collection $(S_m, m \in \mathcal{M})$ of linear subspaces of $\mathbb{H}$. For all $m \in \mathcal{M}$, we can define the estimator $\hat{s}_m$ of $s$, which is the projection estimator of $s$ onto $S_m$. Given some orthonormal basis $(\phi_\lambda, \lambda \in \Lambda_m)$ of $S_m$, it is natural to consider the projection estimator
$$\hat{s}_m = \sum_{\lambda \in \Lambda_m} Y(\phi_\lambda) \phi_\lambda.$$

It is easy to verify that
$$\hat{s}_m = \operatorname{argmin}_{v \in S_m} \left( \|v\|^2 - 2Y(v) \right),$$

which shows that $\hat{s}_m$ does not depend on the particular choice of the basis $(\phi_\lambda, \lambda \in \Lambda_m)$.

It is natural to estimate $T(s)$ by $T(\hat{s}_m)$. Let $s_m$ denote the orthogonal projection of $s$ onto $S_m$. Since $T$ is a linear functional,

$$E\left(T(\hat{s}_m)\right) = T(s_m).$$

Hence, the quadratic risk of the estimator $T(\hat{s}_m)$ can be decomposed into a variance term and a bias term:

$$E\left[(T(\hat{s}_m) - T(s))^2\right] = (T(s_m) - T(s))^2 + E\left[(T(\hat{s}_m) - T(s_m))^2\right].$$

The variance term can be easily computed, by using the properties of the isonormal process $L$.

$$E\left[(T(\hat{s}_m) - T(s_m))^2\right] = \frac{\sigma^2}{n} \sum_{\lambda \in \Lambda_m} T^2(\phi_\lambda) := \sigma_m^2.$$

Our aim is to find an estimator among the collection $(T(\hat{s}_m), m \in \mathcal{M})$ that minimizes the quadratic risk

$$(T(s_m) - T(s))^2 + \sigma_m^2.$$

Model selection by penalized criterion has been introduced by Barron et al. [3] and used in the framework defined by model (1.1) for the estimation of the whole object $s$ by Birgé and Massart [5], and by Laurent and Massart [18] for the estimation of quadratic functionals of $s$. Usually, the bias term, or the sum of this bias term and a term that does not depend on $m \in \mathcal{M}$, is estimated, and the methods proposed in previous papers consist in minimizing over $m \in \mathcal{M}$ this estimation of the bias term, plus some penalty term $\text{pen}(m)$, which has to be suitably chosen. For example, when one estimates $s$, the bias term appearing in the quadratic risk $\mathbb{E}(\|s - \hat{s}_m\|^2)$ equals $\|s - s_m\|^2$. Using Pythagoras'equality, this bias term equals $\|s\|^2 - \|s_m\|^2$. Hence, minimizing $\|s - s_m\|^2 + \text{pen}(m)$ is equivalent to minimize $-\|s_m\|^2 + \text{pen}(m)$, and one can easily find an unbiased estimator of $-\|s_m\|^2$.

In our case, the bias term equals $(T(s_m) - T(s))^2$, and this expression cannot be simplified as previously nor estimated. Therefore, in order to use model selection by penalized criterion methods, we introduce a new criterion. We assume that $\mathcal{M}$ is a subset of $\mathbb{N}$. This implies in particular that $\mathcal{M}$ is ordered. The criterion which is introduced in Definition 1 aims at finding $m \in \mathcal{M}$ which minimizes

$$\sup_{j \geq m, j \in \mathcal{M}} |T(s_m) - T(s_j)| + \sigma_m.$$

**Definition 1** *Let $(S_m, m \in \mathcal{M})$ be a finite or countable collection of linear subspaces of $\mathbb{H}$. For all $m \in \mathcal{M}$, let $(\phi_\lambda, \lambda \in \Lambda_m)$ be an orthonormal basis of $S_m$ and let*

$$\hat{s}_m = \sum_{\lambda \in \Lambda_m} Y(\phi_\lambda) \phi_\lambda.$$

We define, for all $m \in \mathcal{M}$,

$$pen(m) = \sqrt{2x_m}\sigma_m ,$$

where $(x_m, m \in \mathcal{M})$ is a sequence of nonnegative real numbers.

We set, for all $j, m \in \mathcal{M}$,

$$\sigma_{j,m}^2 = \frac{\sigma^2}{n}\mathbb{E}\left[\left(\sum_{\lambda \in \Lambda_m} T(\phi_\lambda)L(\phi_\lambda) - \sum_{\lambda \in \Lambda_j} T(\phi_\lambda)L(\phi_\lambda)\right)^2\right]$$

and

$$H(j,m) = \sqrt{2x_{j,m}}\sigma_{j,m}$$

where $(x_{j,m}, (j,m) \in \mathcal{M}^2)$ is a sequence of nonnegative real numbers.

We define for all $m \in \mathcal{M}$

$$\widehat{Crit}(m) = \sup_{j \geq m, j \in \mathcal{M}} [|T(\hat{s}_m) - T(\hat{s}_j)| - H(j,m)] + pen(m), \quad (2.3)$$

and

$$\hat{m} = \inf\left\{m \in \mathcal{M}, \widehat{Crit}(m) \leq \inf_{j \in \mathcal{M}} \widehat{Crit}(j) + \frac{1}{n}\right\}.$$

We estimate the linear functional $T(s)$ by $T(\hat{s}_{\hat{m}})$.

In the following Theorem, we give an upper bound for the risk with respect to the $\mathbb{L}_p$ loss of the estimator $T(\hat{s}_{\hat{m}})$.

**Theorem 1** *Let $\mathbb{H}$ be some separable Hilbert space endowed with the scalar product $\langle .,. \rangle$. One observes the Gaussian process $\{Y(t), t \in \mathbb{H}\}$, where $Y(t)$ is given by (1.1). Let $T$ be some linear functional defined on $\mathcal{S} \subset \mathbb{H}$. Let $\mathcal{M} \subset \mathbb{N}$ and let $(S_m, m \in \mathcal{M})$ be some finite or countable collection of linear subspaces of $\mathbb{H}$. Let $T(\hat{s}_{\hat{m}})$ be defined in Definition 1. Let for all $m \in \mathcal{M}$,*

$$Crit(m) = \sup_{j \geq m, j \in \mathcal{M}} |T(s_j) - T(s_m)| + pen(m).$$

*Let $m^*$ be defined by*

$$m^* = \inf\left\{m \in \mathcal{M} / Crit(m) \leq \inf_{l \in \mathcal{M}} Crit(l) + \frac{1}{n}\right\}.$$

*Then, for all $p \geq 1$, there exists some positive constant $C(p)$ depending on $p$ only such that*

$$\mathbb{E}\left(|T(\hat{s}_{\hat{m}}) - T(s)|^p\right) \leq C(p)\left((Crit(m^*))^p + |T(s_{m^*}) - T(s)|^p + \sigma_{m^*}^p\right)$$
$$+ C(p)\left(\sup_{j \leq m^*}(H^p(m^*, j)) + \sum_{m \in \mathcal{M}} e^{-x_m}\sigma_m^p + \sum_{j \geq m^*} e^{-x_{j,m^*}}\sigma_{j,m^*}^p + \frac{1}{n^p}\right). \quad (2.4)$$

**Comments:**

- In the definition of $\widehat{\mathrm{Crit}}(m)$ given in (2.3), we compare $T(\hat{s}_m)$ with the estimators $T(\hat{s}_j)$ for $j > m$. This is a common point of our procedure with the initial method from Lepski [19, 20, 21, 23]. An important difference between our procedure and that of Lepski, however, is that we dissociate in (2.3) the terms $H(j,m)$ and $\mathrm{pen}(m)$. This allows us to obtain non asymptotic results based on Gaussian concentration inequalities.
  Let us explain the main ideas underlying the definition of our estimator and how we obtain non asymptotic oracle inequalities for general linear functionals. As mentioned above, our goal is to minimize the criterion

$$\mathrm{Crit}(m) = \sup_{j \geq m, j \in \mathcal{M}} |T(s_j) - T(s_m)| + \mathrm{pen}(m).$$

The first term in this expression is a bias term, and the second one is closely related to the standard deviation of $T(\hat{s}_m)$. The unknown criterion $\mathrm{Crit}(m)$ is estimated by $\widehat{\mathrm{Crit}}(m)$ defined by (2.3). This criterion involves the term $H(j,m)$ which is the standard deviation of $|T(\hat{s}_m) - T(\hat{s}_j)|$ multiplied by $\sqrt{2x_{j,m}}$. For a suitable choice of $x_{j,m}$, we prove that, with high probability,

$$\sup_{j \geq m, j \in \mathcal{M}} [|T(\hat{s}_m) - T(\hat{s}_j)| - H(j,m)] \leq \sup_{j \geq m, j \in \mathcal{M}} |T(s_j) - T(s_m)|.$$

This implies that, with high probability,

$$\forall m \in \mathcal{M}, \widehat{\mathrm{Crit}}(m) \leq \mathrm{Crit}(m).$$

Since $\hat{m}$ is a minimizer of $\widehat{\mathrm{Crit}}(m)$, we obtain that, with high probability,

$$\widehat{\mathrm{Crit}}(\hat{m}) \leq \inf_{m \in \mathcal{M}} \mathrm{Crit}(m),$$

which leads to an oracle inequality. We show that, up to remainder terms of smaller order, the risk of the estimator $T(\hat{s}_{\hat{m}})$ with respect to the $\mathbb{L}_p$ loss behaves as well as the risk of the "best" estimator of the collection. Our procedure is also easily implementable as we will see in the simulation study.

- In order to prove Theorem 1, we do not have to assume that the family $(S_m, m \in \mathcal{M})$ is nested. We just need to order this family. If for example $\mathbb{H} = \mathbb{L}_2([0,1])$, we can mix spaces generated by several kinds of orthonormal bases, for example, a wavelet basis, the Fourier basis, a spline basis. Let us explain how to extend our procedure in the case where we consider $L$ different bases. We set

$$\mathcal{M} = \{(l,m), l \in \mathcal{L}, m \in \mathcal{M}_l\},$$

where $\mathcal{L} = \{1, 2, \ldots, L\}$ and $\mathcal{M}_l \subset \mathbb{N}$. For all $m \in \mathcal{M}_l$, and $l \in \{1, 2, \ldots, L\}$, we set

$$\widehat{\mathrm{Crit}}_l(m) = \sup_{j \geq m, j \in \mathcal{M}_l} [|T(\hat{s}_m) - T(\hat{s}_j)| - H(j,m)] + \mathrm{pen}(m).$$

We define

$$(\hat{l}, \hat{m}) = \inf \left\{ (l, m) \in \mathcal{M}, \widehat{\mathrm{Crit}}_l(m) \leq \inf_{(k,j) \in \mathcal{M}} \widehat{\mathrm{Crit}}_k(j) + \frac{1}{n} \right\},$$

where $\mathcal{M}$ is ordered by the lexicographical order.
- It would be simpler to define $\hat{m}$ as a minimizer of $\widehat{\mathrm{Crit}}(m)$, but such a minimizer may not exist, or not be unique. This explains why we add the term $1/n$ in the definition of $\hat{m}$ given in Definition 1.

We shall derive in the next section applications of Theorem 1 to adaptive results in the minimax sense for pointwise estimation and for the estimation of the mean of the signal on some interval.

## 3. Minimax results

### 3.1. Pointwise adaptive estimation

Assume we observe $(Y(u), u \in [0, 1])$ which obeys the Gaussian white noise model:

$$Y(u) = \int_0^u s(x)dx + \frac{\sigma}{\sqrt{n}} W(u), u \in [0, 1], \quad (3.5)$$

where $s \in \mathbb{H} = \mathbb{L}_2([0,1])$ and $W$ is a standard Brownian motion.

Let $r \geq 0$, we assume that $s^{(r)}$ exists and that $s^{(r)} \in \mathcal{C}([0, 1])$, the set of continuous functions on $[0, 1]$. We consider the problem of estimating $T(s) = s^{(r)}(x_0)$ for some fixed $x_0 \in [0, 1]$.

We introduce the following notation: let $\{S_j, j \geq 0\}$ be a multiresolution analysis with father wavelet $\varphi$ and mother wavelet $\psi$ (see for example [15]). Define

$$\begin{aligned}
\varphi_{j,k}(x) &= 2^{j/2} \varphi(2^j x - k), \, x \in [0, 1], \, j \geq 0 \text{ and } k \in \mathbb{Z}; \\
\psi_{j,k}(x) &= 2^{j/2} \psi(2^j x - k), \, x \in [0, 1], \, j \geq 0 \text{ and } k \in \mathbb{Z}.
\end{aligned}$$

For all $m \geq 0$, $S_m$ denotes the linear space spanned by the functions $(\varphi_{m,k}, k \in \mathbb{Z})$. We recall that $s_m$ denotes the orthogonal projection of $s$ onto $S_m$:

$$s_m = \sum_{k \in \mathbb{Z}} \langle s, \varphi_{m,k} \rangle \varphi_{m,k},$$

and that $\hat{s}_m$ denotes the estimator of $s$ based on the model $S_m$:

$$\hat{s}_m = \sum_{k \in \mathbb{Z}} \int_0^1 \varphi_{m,k}(u) dY(u) \varphi_{m,k}.$$

We will consider compactly supported wavelets $\varphi$ for which the above sum is finite.

For all $\alpha > 0$ and $1 \leq q \leq +\infty$, the notation $([0,1])$ is used for the classical Besov space endowed with the norm $\|\cdot\|_{\alpha,\infty,q}$ (see for example [15, Definition 9.2]). We denote by $B_{\alpha,\infty,q}(L)$ the set of functions $s$ in $B_{\alpha,\infty,q}([0,1])$ such that $\|s\|_{\alpha,\infty,q} \leq L$.

Assume that the following conditions are satisfied by $\varphi$ and $\psi$:

**(i)** $\exists M \geq 0$ such that $\text{supp}(\varphi)$ and $\text{supp}(\psi)$ are included in $[-M, M]$.
**(ii)** $\exists K \geq 0$ such that $\|\varphi\|_\infty \vee \|\psi\|_\infty \leq K$.
**(iii)** $\exists N \geq 0$ such that $\int x^n \psi(x) dx = 0$ for $n = 0, \ldots, N$.
**(iv)** We assume that $\varphi^{(r)}$ exists and is bounded on $\text{supp}(\varphi)$ by $K_r$.

**Corollary 1** *Let $d_n$ denote the integer part of $\ln(n)/\ln(2)$ and let $\mathcal{M} = \{1, \ldots, d_n\}$. For all $m \in \mathcal{M}$, let $S_m$ be the linear space spanned by the functions $(\varphi_{m,k}, k \in \mathbb{Z})$.*

*For $p \geq 1$, we define*

$$\forall m \in \mathcal{M}, x_m = \frac{p}{2} \ln\left(2^{m(1+2r)}\right),$$

$$\forall (j,m) \in \mathcal{M}^2 \text{ if } j > m, x_{j,m} = \frac{p}{2} \ln\left(2^{j(1+2r)} - 2^{m(1+2r)}\right), \quad x_{m,m} = 0.$$

*Let $\hat{m}$ be defined as in Definition 1. Let $r < \alpha$.*

*There exist some constants $C$ depending on $\alpha$, $p$, $q$ and $r$ and $C(\sigma)$ depending on $\sigma$ such that, for any integer $n$ satisfying $nL^2/(\sigma^2 \ln(n)) \geq 2^{1+2\alpha}$ and $n^{2\alpha}\sigma^2 \ln(n) \geq L^2$,*

$$\sup_{s \in B_{\alpha,\infty,q}(L)} \mathbb{E}\left(\left|\hat{s}_{\hat{m}}^{(r)}(x_0) - s^{(r)}(x_0)\right|^p\right) \leq CL^{\frac{p(1+2r)}{1+2\alpha}} \left(\frac{\sigma^2 \ln n}{n}\right)^{\frac{p(\alpha-r)}{2\alpha+1}}$$

$$+ C(\sigma) \frac{\ln n}{n^{p/2}}.$$

Comments on the optimality of the result stated in Corollary 1 are given in Subsection 3.3, where the multidimensional case is considered.

### 3.2. Estimation of the mean of the signal on an interval

As above, we observe $(Y(u), u \in [0,1])$ defined by (3.5). We use the same notation as in Section 3.1. We now consider the problem of estimating the linear functional

$$T(s) = \frac{1}{H_n} \int_{I_{H_n}} s(x) dx$$

where $I_{H_n}$ is an interval included in $[0,1]$ with length $H_n$ that may depend on $n$.

**Corollary 2** *Let $(Y(u), u \in [0,1])$ defined by (3.5). Let $I_{H_n}$ be some interval included in $[0,1]$ with length $H_n$. Let $T(s) = \int_{I_{H_n}} s(x) dx/H_n$. Let $m_n = \sup\{m \in \mathbb{N}, 2^m \leq 1/H_n\}$. Let $\mathcal{M} = \{0, 1, \ldots m_n + 1\}$. For all $m \in \mathcal{M}\setminus\{m_n + 1\}$, let $S_m$ be the linear space spanned by the functions $(\varphi_{m,k}, k \in \mathbb{Z})$ and let $S_{m_n+1}$ be the linear space spanned by the indicator function of the interval $I_{H_n}$.*

For $p \geq 1$, we define

$$x_m = \frac{pm}{2} \quad \forall m \in \mathcal{M}, m \leq m_n, \ x_{m_n+1} = \frac{p}{2}\ln(1/H_n)$$

$$x_{j,m} = \frac{p(j \vee m)}{2} \quad \forall (j,m) \in \mathcal{M}^2, j \neq m \ \text{and} \ x_{m,m} = 0$$

Let $\hat{m}$ be defined as in Definition 1. Let $n$ be any integer satisfying $nL^2/\sigma^2 \ln(n) \geq 1$.

Then, for all $\alpha > 0, q \geq 1, p \geq 1$, there exist some constants $C$ depending on $\alpha$, $p$ and $q$ and $C(\sigma)$ depending on $\sigma$ such that the following inequalities hold:

If $H_n \leq \left(\sigma^2 \ln(n)/nL^2\right)^{1/(1+2\alpha)}$,

$$\sup_{s \in B_{\alpha,\infty,q}(L)} \mathbb{E}\left(|T(\hat{s}_{\hat{m}}) - T(s)|^p\right) \leq CL^{\frac{p}{1+2\alpha}} \left(\frac{\sigma^2 \ln n}{n}\right)^{\frac{p\alpha}{2\alpha+1}} + \frac{C(\sigma)}{n^{p/2}}.$$

If $H_n \geq \left(\sigma^2 \ln(n)/nL^2\right)^{1/(1+2\alpha)}$,

$$\sup_{s \in B_{\alpha,\infty,q}(L)} \mathbb{E}\left(|T(\hat{s}_{\hat{m}}) - T(s)|^p\right) \leq C\frac{\sigma^p}{(nH_n)^{p/2}}\left(\ln(1/H_n)\right)^{p/2} + \frac{C(\sigma)}{n^{p/2}}.$$

**Comments:**

- The rates of convergence that we obtain depend on the relation between $H_n$, $n$, $\sigma$ and the regularity of the signal (via the parameters $\alpha$ and $L$). If $H_n \geq \left(\sigma^2\ln(n)/nL^2\right)^{1/(1+2\alpha)}$, the best estimator is the "naif" estimator $\int_{I_{H_n}} dY(u)/H_n$, which is unbiased and achieves the rate $1/\sqrt{nH_n}$. Note that in our result, we loose a logarithmic term $(\ln(1/H_n))$ for the adaptation to the unknown regularity of the signal. When $H_n$ is independent of $n$, we recover the parametric rate $1/\sqrt{n}$ for the estimation of $T(s)$.
  If $H_n \leq \left(\sigma^2\ln(n)/nL^2\right)^{1/(1+2\alpha)}$, we obtain the same rates for the estimation of $T(s)$ as for the estimation of the signal $s$ at one point. In this case, the "naif" estimator has a too large variance, and one takes advantage of considering estimators which are biased, but with smaller variance.
- Our procedure is adaptive with respect to the unknown link between $H_n$ and the regularity of the signal and allows us to obtain the optimal rates (up to logarithmic terms due to adaptation) in both cases as explained below.
- It was possible to establish the upper bounds given in Corollary 2 since the result stated in Theorem 1 is non asymptotic. We have indeed considered here a linear functional that depends on $n$.

**Lower bounds:**

Using the results of Donoho and Liu [13], one can show that, up to logarithmic terms, the upper bounds given in Corollary 2 are optimal over Hölderian balls. The lower bounds are given in the following lemma.

**Lemma 1** *Let $(Y(u), u \in [0,1])$ defined by (3.5). Let $0 \leq H_n \leq 1/2$ and $I_{H_n} = [0, H_n]$, we consider the linear functional $T(s) = \int_{I_{H_n}} s(x)dx/H_n$. Set, for all $\alpha \in ]0,1]$ and $L > 0$,*

$$\mathcal{H}_\alpha(L) = \{f : [0,1] \to \mathbb{R}, \; \forall x, y \in [0,1], \; |f(x) - f(y)| \leq L|x-y|^\alpha\}.$$

*If $H_n \geq \big((1+2\alpha)/(2nL^2)\big)^{1/(1+2\alpha)}$,*

$$\inf_{T_n} \sup_{s \in \mathcal{H}_\alpha(L)} \mathbb{E}\left((T_n - T(s))^2\right) \geq \frac{C(\alpha, \sigma)}{nH_n} \tag{3.6}$$

*where $C(\alpha, \sigma)$ is a constant depending on $\alpha$ and $\sigma$ and where the infimum is taken over all possible estimators.*

*If $H_n \leq \big((1+2\alpha)/(nL^2)\big)^{1/(1+2\alpha)}/2$, and $nL^2 \geq 1 + 2\alpha$,*

$$\inf_{T_n} \sup_{s \in \mathcal{H}_\alpha(L)} \mathbb{E}\left((T_n - T(s))^2\right) \geq C(\alpha, \sigma) L^{2/(1+2\alpha)} n^{-2\alpha/(1+2\alpha)}. \tag{3.7}$$

### 3.3. Multidimensional pointwise adaptive estimation

One observes the Gaussian white noise model

$$Z(x) = \int_{[0,1]^d} \mathbb{1}_{[0,x_1] \times \cdots \times [0,x_d]}(u)s(u)du + \frac{1}{\sqrt{n}}W(x)$$

for all $x = (x_1, \ldots, x_d) \in [0,1]^d$, where $W$ is the standard Wiener Process on $[0,1]^d$.

We consider $\mathbb{H} = \mathbb{L}_2([0,1]^d)$ endowed with its usual scalar product. We assume that $s \in \mathcal{C}([0,1]^d)$, the set of continuous functions on $[0,1]^d$. Let $x_0 \in [0,1]^d$, we estimate $T(s) = s(x_0)$.

Our aim is to obtain adaptive results in the minimax sense over isotropic Hölder spaces defined as follows. For all $\alpha \in ]0,1]$ and $L > 0$, let

$$\mathcal{H}_\alpha(L) = \{f : [0,1]^d \to \mathbb{R}, \forall x, y \in [0,1]^d, |s(x) - s(y)| \leq L\|x-y\|_\infty^\alpha\},$$

where $\|x - y\|_\infty = \sup_{1 \leq i \leq n} |x_i - y_i|$. In order to estimate $T(s)$, we use the Haar basis of $\mathbb{L}_2([0,1]^d)$.

For all $m \in \mathcal{M}$, let $S_m$ be the space of piecewise constant functions on the sets $[\frac{k_1}{2^m}, \frac{k_1+1}{2^m}[ \times \cdots \times [\frac{k_d}{2^m}, \frac{k_d+1}{2^m}[$ for all $(k_1, \ldots, k_d) \in \{0, 1, \ldots, 2^m - 1\}^d$.

**Corollary 3** *Let $d_n$ denote the integer part of $\ln(n)/(d \ln(2))$ and let $\mathcal{M} = \{1, \ldots, d_n\}$.*

*For $p \geq 1$, we define*

$$\forall m \in \mathcal{M}, x_m = \frac{pd}{2} \ln(2^m)$$

$$\forall (j, m) \in \mathcal{M}^2, x_{j,m} = \frac{pd}{2} \ln(2^j).$$

Let $\hat{m}$ be defined in Theorem 1. For all $\alpha \in ]0,1]$ and $L > 0$, the following inequality holds if $n/(\sigma^2 \ln(n)) \geq 2^{d+2\alpha} L^{-2}$ and $L^2 \leq n^{\frac{2\alpha}{d}} \sigma^2 \ln(n)$ :

$$\sup_{s \in \mathcal{H}_\alpha(L)} \mathbb{E}\left(|\hat{s}_{\hat{m}}(x_0) - s(x_0)|^p\right) \leq C L^{\frac{pd}{d+2\alpha}} \left(\frac{\sigma^2 \ln n}{n}\right)^{\frac{p\alpha}{2\alpha+d}} + C(\sigma)\frac{\ln n}{n^{p/2}},$$

where $C$ is a constant depending on $p, \alpha$ and $d$ and $C(\sigma)$ is a constant depending on $\sigma$.

**Comments.** It follows from the results given in Lepski [19] and Brown and Low [7] that the rates obtained in Corollary 1 and in Corollary 3 are optimal. These authors showed that the logarithmic loss which appears in the rate of convergence compared with the minimax rates is unavoidable for adaptive estimators.

The dependence of our upper bound for the risk with respect to the radius $L$ of the Besov or Hölderian balls is the sharp one obtained by Klemelä and Tsybakov [17].

## 4. Simulation study

Throughout this section, we consider the finite dimensional Gaussian regression model. The regression functions that we consider are defined on $[0,1]$ by:

$$s_1(x) = (x^4 - x)\sin(6x),$$

$$s_2(x) = \exp(-30|x - 0.75|) + \exp(-30|x - 0.25|),$$

$$s_3(x) = x\cos(2\pi x)\; \mathbb{1}_{0 < x \leq 2/3} + x^2 \cos(15\pi x)\; \mathbb{1}_{2/3 < x \leq 1}.$$

The estimation is based on the simulations

$$y_i = s_j\left(\frac{i}{n}\right) + \sigma \varepsilon_i \quad i = 1, \ldots, n \quad j = 1, 2, 3 \qquad (4.8)$$

with $(\varepsilon_1, \ldots, \varepsilon_n)$ independent standard normal variables, $\sigma = 0.2$ and $n = 256$.

We set $d_n = \ln_2(n) = 8$ and $\mathcal{M} = \{1, 2, \ldots, d_n\}$. The estimators are built using a wavelet basis. We use the Haar basis (denoted by H in the tables) or the Daubechies 20 basis (denoted by D 20). In both cases, for all $m \in \mathcal{M}$, $S_m$ is the linear space spanned by the functions $(\varphi_{m,k}, k \in \mathbb{Z})$, where $\varphi_{m,k} = 2^{m/2}\varphi(2^m \cdot - k)$, $\varphi$ is the father wavelet of the basis. The results presented in the tables must be divided by 100. All simulations were programmed in Matlab 7.3 with the wavelab wavelet toolbox.

### 4.1. Pointwise estimation

We first consider the estimation of the linear functional $T(s) = s(\tilde{x}_l)$ for some fixed points $\tilde{x}_l \in [0,1]$.

When the Haar basis is used to construct the estimators, we obtain

$$\forall m \in \mathcal{M}, \sigma_m^2 = 2^m \frac{\sigma^2}{n}, \quad \forall\, 1 \leq m \leq j \leq d_n, \sigma_{j,m}^2 = (2^j - 2^m)\frac{\sigma^2}{n}.$$

We set

$$x_{j,m} = \frac{1}{2}\ln(2^j - 2^m),$$

$$H(j,m) = (2^j - 2^m)^{1/2}\sqrt{2x_{j,m}}\frac{\sigma}{\sqrt{n}},$$

and

$$x_m = \frac{1}{2}\ln(2^m),$$

$$\text{pen}(m) = \sqrt{2x_m}\,2^{m/2}\frac{\sigma}{\sqrt{n}}.$$

The choices of $H(j,m)$ and of $\text{pen}(m)$ given above correspond to the control of the $\mathbb{L}_1$ risk in Corollary 1 ($p=1$), when we use the Haar basis. For these choices,

$$\sum_{m \in \mathcal{M}} e^{-x_m}\sigma_m \leq \frac{\sigma \ln_2(n)}{\sqrt{n}}$$

and for all $m \in \mathcal{M}$,

$$\sum_{j \geq m} e^{-x_{j,m}}\sigma_{j,m} \leq \frac{\sigma \ln_2(n)}{\sqrt{n}}.$$

The order of magnitude of both series is smaller than the rates of convergence of $\mathbb{E}(|\hat{s}(x) - s(x)|)$ obtained in Corollary 1.

Our procedure is called **P1**. We compare the performances of our procedure to the performances of the estimator $\tilde{s}$ studied by Baraud [1] and defined as follows: let, for all functions $t$

$$\gamma_n(t) = \frac{1}{n}\sum_{i=1}^{n}\left(y_i - t\left(\frac{i}{n}\right)\right)^2,$$

$$\tilde{s} = \operatorname{argmin}_{m \in \mathcal{M}}(\gamma_n(\hat{s}_m) + \text{pen}'(m)),$$

with $\text{pen}'(m) = 2.2^m \sigma^2/n$, which corresponds to a Mallow's $C_p$ criterion. This procedure is called **P2**.

We also compare our procedure with a wavelet thresholding procedure, for which the wavelet coefficients which are smaller than $\sigma\sqrt{2\ln(n)}$ are set to 0. This procedure is called **P3**.

In Figure 1, we have represented the functions $s_1, s_2, s_3$ and one simulated sample for the noised observations $(i/n, y_i)$.

We estimate the pointwise risk in absolute value $\mathbb{E}(|\hat{s}(x) - s(x)|)$ for the estimation of $s_j(x)$ with the procedures **P1, P2** and **P3**. The estimation of the risk is based on $N = 5000$ simulations and is defined as

$$\hat{r}_1(x) = \frac{1}{N}\sum_{l=1}^{N}\left|\hat{s}^{(l)}(x) - s(x)\right|$$

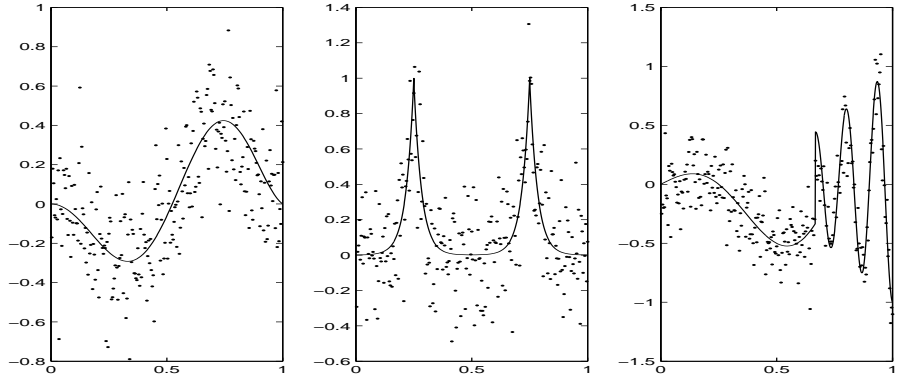

FIG 1. *Functions $s_1, s_2, s_3$, and one simulated sample $n = 2^8$, $\sigma = 0.2$.*

where $s = s_1, s_2$ or $s_3$ and $\hat{s}^{(l)}$ is the estimator of $s$ based on the $l$-th simulated sample. In the following tables, we give the values of $100 * \hat{r}_1$ at points $1/4$, $1/3$, $1/2$ and $3/4$ for $s_1$, $1/8$, $1/4$, $1/3$ and $1/2$ for $s_2$ and $1/4$, $1/3$, $1/2$ and $7/8$ for $s_3$.

| $s_1$ | $\tilde{x}_1 = 1/4$ | | | $\tilde{x}_2 = 1/3$ | | | $\tilde{x}_3 = 1/2$ | | | $\tilde{x}_4 = 3/4$ | | |
|---|---|---|---|---|---|---|---|---|---|---|---|---|
| | **P1** | **P2** | **P3** | **P1** | **P2** | **P3** | **P1** | **P2** | **P3** | **P1** | **P2** | **P3** |
| H | 5.6 | 3.2 | 14.9 | 4.5 | 4.0 | 7.4 | 4.2 | 6.9 | 11.3 | 5.7 | 8.0 | 16.9 |
| D20 | 5.4 | 5.4 | 5.8 | 2.6 | 3.0 | 3.0 | 6.5 | 6.5 | 2.0 | 6.6 | 6.6 | 7.5 |

| $s_2$ | $\tilde{x}_1 = 1/8$ | | | $\tilde{x}_2 = 1/4$ | | | $\tilde{x}_3 = 1/3$ | | | $\tilde{x}_4 = 1/2$ | | |
|---|---|---|---|---|---|---|---|---|---|---|---|---|
| | **P1** | **P2** | **P3** | **P1** | **P2** | **P3** | **P1** | **P2** | **P3** | **P1** | **P2** | **P3** |
| H | 3.8 | 6.3 | 3.2 | 23.3 | 27.8 | 30.4 | 4.7 | 6.3 | 4.8 | 3.5 | 6.1 | 3.0 |
| D20 | 6.8 | 5.1 | 5.6 | 20.8 | 23.3 | 35.8 | 6.1 | 6.5 | 9.7 | 6.8 | 5.0 | 6.7 |

| $s_3$ | $\tilde{x}_1 = 1/4$ | | | $\tilde{x}_2 = 1/3$ | | | $\tilde{x}_3 = 1/2$ | | | $\tilde{x}_4 = 7/8$ | | |
|---|---|---|---|---|---|---|---|---|---|---|---|---|
| | **P1** | **P2** | **P3** | **P1** | **P2** | **P3** | **P1** | **P2** | **P3** | **P1** | **P2** | **P3** |
| H | 5.9 | 7.9 | 5.9 | 5.2 | 8.0 | 5.0 | 8.0 | 7.9 | 9.9 | 7.5 | 8.2 | 8.1 |
| D20 | 3.4 | 5.2 | 4.9 | 4.9 | 6.3 | 6.1 | 3.6 | 5.1 | 6.3 | 5.2 | 5.2 | 7.1 |

For our procedure, at each point $\tilde{x}_l$, an estimator is selected among a collection of $d_n$ estimators. This collection is composed of the estimators based on a projection onto a wavelet basis up to the level $j$ for $j = 1, \ldots, d_n$. We represent in Figure 2 the histograms of the selected levels for the estimation of $s_2(x)$ with the Haar basis, for a point where the function $s_2$ is nearly flat ($\tilde{x}_4 = 1/2$) and at a peak of the function $s_2$ ($\tilde{x}_2 = 1/4$). We also represent the histogram of the selected levels for the estimation of $s_2$ with the procedure **P2** when we use the Haar basis. We recall that for this procedure a level is selected for the estimation of the whole function.

Figure 2 clearly shows that, as expected, the level which is selected by our procedure is higher at points where the function to be estimated is "irregular".

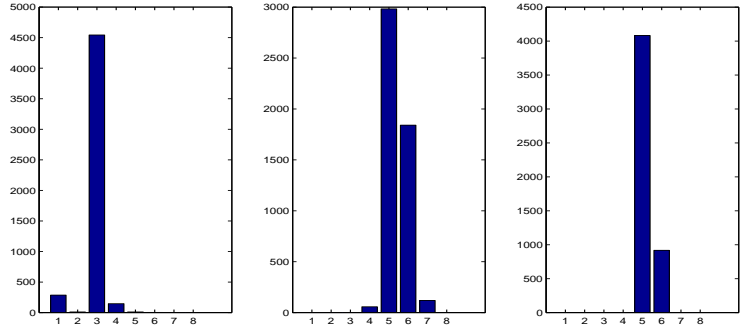

FIG 2. *Histogram of the selected levels for the procedure P1 and $s_2$ at $\tilde{x}_4 = 1/2$ (on the left) and $\tilde{x}_2 = 1/4$ (on the middle), selected levels for the procedure P2 (on the right).*

The procedure **P2** selects a high level to estimate accurately the function $s_2$ near the peaks.

The simulation results show that, in most cases, we obtain good results with the pointwise adaptive procedure **P1** for the risk $\hat{r}_1(x)$. Our procedure performs better at points where the function is "irregular". Except for the fonction $s_1$ with the Haar basis, the procedure **P1** performs in most cases better that **P2**. Whatever the basis and the function, our procedure performs better in most cases than the procedure **P3**, or as well as **P3**.

## 4.2. Estimation of integral functionals

One observes $(y_i, 1 \leq i \leq n)$ given by (4.8). We consider the problem of estimating $\int_0^1 s_j(x)g(x)dx$. In the first part of the study, $g = \mathbb{1}_{[0,H]}/H$, where $H = 1/4, 1/32, 1/128$. The last value of $H$ is comparable to $1/n$. This problem has been considered in Section 3.2.

In the second part of the study, $g$ equals $g_1$ or $g_2$ defined on $[0,1]$ by

$$g_1(x) = \cos(64\pi x)$$
$$g_2(x) = \cos(4\pi x).$$

For all $(j, m) \in \mathcal{M}^2$, we choose the same values for $x_m$ and $x_{j,m}$ as in Section 4.1. Denoting, for all $m \in \mathcal{M}$, by $\pi_{S_m}$ the orthogonal projection onto $S_m$ one has

$$\sigma_m^2 = \|\pi_{S_m}(g)\|^2 \quad \forall m \in \mathcal{M},$$
$$\sigma_{j,m}^2 = \|\pi_{S_j}(g)\|^2 - \|\pi_{S_m}(g)\|^2 \quad \forall j \geq m \in \mathcal{M}.$$

We compare the estimator obtained with our procedure **P1** with the estimators $\int_0^1 g\tilde{s}$ where $\tilde{s}$ is obtained by Procedure **P2** or **P3**. Those two procedures are still denoted by **P2** and **P3**. We also compute the empirical estimator $\sum_{i=1}^{n} y_i g(x_i)/n$. This procedure is denoted by **P4**. Throughout the section, we use the Haar basis for the simulations. In the following tables, we give the value of $100 * \hat{r}_1$,

$$\hat{r}_1 = \frac{1}{N} \sum_{l=1}^{N} \left| \hat{T}(s)^{(l)} - T(s) \right|$$

with $N = 5000$ and $\hat{T}(s)^{(l)}$ the estimation of $T(s)$ based on the $l$-th simulated sample.

|       | $H = 1/4$ |     |     |     | $H = 1/32$ |     |     |     |
|-------|-----|-----|-----|-----|-----|-----|-----|-----|
|       | **P1** | **P2** | **P3** | **P4** | **P1** | **P2** | **P3** | **P4** |
| $s_1$ | 2.1 | 2.0 | 3.5 | 2.0 | 4.6 | 4.2 | 9.3 | 5.7 |
| $s_2$ | 1.7 | 2.1 | 1.4 | 2.1 | 3.3 | 5.7 | 2.6 | 5.7 |
| $s_3$ | 2.9 | 2.0 | 2.0 | 2.0 | 4.4 | 5.7 | 4.4 | 5.7 |

|       | $H = 1/128$ |     |     |     |
|-------|-----|-----|-----|-----|
|       | **P1** | **P2** | **P3** | **P4** |
| $s_1$ | 4.7 | 4.1 | 9.4 | 11.2 |
| $s_2$ | 3.4 | 6.2 | 2.6 | 11.3 |
| $s_3$ | 5.3 | 8.1 | 5.4 | 11.2 |

In most cases, our procedure is comparable with the best procedure. Whatever the value of $H$, the procedures **P2** and **P3** use the same estimator for the function $s$. The risk for **P4** increases as $H$ becomes smaller since this procedure considers the mean of the observations over a smaller sample. Our procedure **P1** takes advantage of the regularity of the signal in a neighbourhood of the interval $[0, H]$ to consider the mean over a larger sample.

In the following tables, we present the results for the estimation of the linear functionals $\int_0^1 s_j(x) g_i(x) dx$, $i = 1, 2$.

|       | $g_1$ |     |     |     |
|-------|-----|-----|-----|-----|
|       | **P1** | **P2** | **P3** | **P4** |
| $s_1$ | $2.09 * 10^{-3}$ | $2.24 * 10^{-3}$ | $2.71 * 10^{-2}$ | 0.7 |
| $s_2$ | 0.29 | 0.30 | 0.28 | 0.72 |
| $s_3$ | 0.36 | 0.30 | 0.31 | 0.74 |
|       | $g_2$ |     |     |     |
|       | **P1** | **P2** | **P3** | **P4** |
| $s_1$ | 2.86 | 2.86 | 2.86 | 2.88 |
| $s_2$ | 0.77 | 0.72 | 1.00 | 0.72 |
| $s_3$ | 0.71 | 0.71 | 0.56 | 0.72 |

Our procedure is comparable to **P2** and in most cases our risk has the same order as that of the best procedure.

## 5. Proofs

### 5.1. Proof of Theorem 1

We shall use the following lemma.

**Lemma 2** *For all $m \in \mathcal{M}$, for all $x > 0$,*

$$\mathbb{P}\left(\widehat{Crit}(m) > Crit(m) + \sqrt{2x}\right) \leq \sum_{j \geq m, j \in \mathcal{M}} e^{-x_{j,m}} e^{-x/\sigma_{j,m}^2}.$$

#### 5.1.1. Proof of Lemma 2

We recall that, for all $m \in \mathcal{M}$,

$$\hat{s}_m = \sum_{\lambda \in \Lambda_m} Y(\phi_\lambda) \phi_\lambda.$$

Since $T$ is a linear functional,

$$T(\hat{s}_m) = \sum_{\lambda \in \Lambda_m} Y(\phi_\lambda) T(\phi_\lambda).$$

Moreover,

$$Y(\phi_\lambda) \sim \mathcal{N}(\langle s, \phi_\lambda \rangle, \frac{\sigma^2}{n}).$$

Using the properties of the isonormal process $L$,

$$T(\hat{s}_j) - T(\hat{s}_m) \sim \mathcal{N}(T(s_j) - T(s_m), \sigma_{j,m}^2).$$

Let $X \sim \mathcal{N}(\mu, v^2)$. For all $x > 0$,

$$\mathbb{P}\left(|X - \mu| > \sqrt{2x}\right) \leq e^{-x/v^2}$$

which implies that

$$\mathbb{P}\left(|X| > |\mu| + \sqrt{2x}\right) \leq e^{-x/v^2}. \tag{5.9}$$

Since for all $a, b > 0$, $\sqrt{a+b} \leq \sqrt{a} + \sqrt{b}$, we obtain that for all $(j,m) \in \mathcal{M}^2$ such that $j \geq m$, for all $x > 0$,

$$\mathbb{P}\left(|T(\hat{s}_j) - T(\hat{s}_m)| \geq |T(s_j) - T(s_m)| + \sqrt{2x_{j,m}} \sigma_{j,m} + \sqrt{2x}\right)$$
$$\leq \mathbb{P}\left(|T(\hat{s}_j) - T(\hat{s}_m)| \geq |T(s_j) - T(s_m)| + \sigma_{j,m}\sqrt{2x_{j,m} + 2x/\sigma_{j,m}^2}\right)$$
$$\leq e^{-x_{j,m}} e^{-x/\sigma_{j,m}^2}.$$

Finally,

$$\mathbb{P}\left(\widehat{\mathrm{Crit}}(m) > \mathrm{Crit(m)} + \sqrt{2x}\right)$$
$$\leq \mathbb{P}\left(\exists j \geq m, j \in \mathcal{M}, |T(\hat{s}_j) - T(\hat{s}_m)| - H(j,m) \geq |T(s_j) - T(s_m)| + \sqrt{2x}\right)$$
$$\leq \sum_{j \geq m, j \in \mathcal{M}} \mathbb{P}\left(|T(\hat{s}_j) - T(\hat{s}_m)| - H(j,m) \geq |T(s_j) - T(s_m)| + \sqrt{2x}\right).$$

This concludes the proof of Lemma 2.

• We first consider the case where $\hat{m} \leq m^*$. Since for all $m \in \mathcal{M}$, $\mathrm{pen}(m) \geq 0$ and by definition of $\widehat{\mathrm{Crit}}(m)$, we have

$$\begin{aligned}\widehat{\mathrm{Crit}}(\hat{m}) &\geq |T(\hat{s}_{\hat{m}}) - T(\hat{s}_{m^*})| - H(m^*, \hat{m}) \\ &\geq |T(\hat{s}_{\hat{m}}) - T(s)| - |T(\hat{s}_{m^*}) - T(s)| - H(m^*, \hat{m}).\end{aligned}$$

Since

$$\widehat{\mathrm{Crit}}(\hat{m}) \leq \widehat{\mathrm{Crit}}(m^*) + \frac{1}{n}$$

we obtain

$$|T(\hat{s}_{\hat{m}}) - T(s)| \leq \widehat{\mathrm{Crit}}(m^*) + H(m^*, \hat{m}) + |T(\hat{s}_{m^*}) - T(s)| + \frac{1}{n}.$$

On the event $\{\hat{m} \leq m^*\}$,

$$H(m^*, \hat{m}) \leq \sup_{j \leq m^*} H(m^*, j).$$

Hence, using Lemma 2, we obtain that for all $x > 0$, the probability of

$$\left\{|T(\hat{s}_{\hat{m}}) - T(s)| > \mathrm{Crit}(m^*) + \sqrt{2x} + \sup_{j \leq m^*} H(m^*, j) + |T(\hat{s}_{m^*}) - T(s)| + \frac{1}{n}\right\}$$

$\cap \{\hat{m} \leq m^*\}$ is bounded by

$$\sum_{j \geq m^*} e^{-x_{j,m^*}} e^{-x/\sigma^2_{j,m^*}}. \tag{5.10}$$

• We now consider the case where $\hat{m} > m^*$. We recall that

$$\sigma_m^2 = \mathrm{var}(T(\hat{s}_m)) = \frac{\sigma^2}{n} \sum_{\lambda \in \Lambda_m} T^2(\phi_\lambda).$$

Using inequality (5.9), we obtain that for all $m \in \mathcal{M}$,

$$\mathbb{P}\left(|T(\hat{s}_m) - T(s)| \geq |T(s_m) - T(s)| + \sqrt{2x} + \sigma_m \sqrt{2x_m}\right)$$
$$\leq e^{-x_m} e^{-x/\sigma_m^2}.$$

This implies that

$$\mathbb{P}\left(|T(\hat{s}_{\hat{m}}) - T(s)| \geq |T(s_{\hat{m}}) - T(s)| + \sqrt{2x} + \mathrm{pen}(\hat{m})\right)$$
$$\leq \sum_{m \in \mathcal{M}} e^{-x_m} e^{-x/\sigma_m^2}.$$

We notice that for all $m \in \mathcal{M}$,

$$\widehat{\mathrm{Crit}}(m) \geq \mathrm{pen}(m)$$

since

$$\sup_{j \geq m} [|T(\hat{s}_m) - T(\hat{s}_j)| - H(j, m)] \geq [|T(\hat{s}_m) - T(\hat{s}_m)| - H(m, m)]$$

and the right hand side of the inequality is equal to 0.

Using the inequalities

$$\mathrm{pen}(\hat{m}) \leq \widehat{\mathrm{Crit}}(\hat{m}) \leq \widehat{\mathrm{Crit}}(m^*) + \frac{1}{n},$$

we obtain

$$\mathbb{P}\left(|T(\hat{s}_{\hat{m}}) - T(s)| \geq |T(s_{\hat{m}}) - T(s)| + \sqrt{2x} + \widehat{\mathrm{Crit}}(m^*) + \frac{1}{n}\right)$$
$$\leq \sum_{m \in \mathcal{M}} e^{-x_m} e^{-x/\sigma_m^2}.$$

We use the following inequality that holds if $\hat{m} > m^*$:

$$|T(s_{\hat{m}}) - T(s)| \leq \sup_{j \geq m^*} |T(s_j) - T(s)|$$

and we apply Lemma 2 with $m = m^*$ to control $\widehat{\mathrm{Crit}}(m^*)$ by $\mathrm{Crit}(m^*)$. Hence the probability of

$$\left\{|T(\hat{s}_{\hat{m}}) - T(s)| \geq \sup_{j \geq m^*} |T(s_j) - T(s)| + 2\sqrt{2x} + \mathrm{Crit}(m^*) + \frac{1}{n}\right\}$$
$$\cap \{\hat{m} > m^*\} \text{ is bounded by}$$

$$\sum_{m \in \mathcal{M}} e^{-x_m} e^{-x/\sigma_m^2} + \sum_{j \geq m^*, j \in \mathcal{M}} e^{-x_{j,m^*}} e^{-x/\sigma_{j,m^*}^2} . \quad (5.11)$$

Define

$$C_{m^*} = \mathrm{Crit}(m^*) + \sup_{j \leq m^*} H(m^*, j) + \sup_{j \geq m^*} |T(s_j) - T(s)| + \frac{1}{n}$$

and

$$X = |T(\hat{s}_{\hat{m}}) - T(s)| , \; Y = |T(\hat{s}_{m^*}) - T(s)| .$$

It follows from (5.10) and (5.11) that for all $x > 0$,

$$\mathbb{P}\left(X - Y > C_{m^*} + 2\sqrt{2x}\right) \leq \sum_{m \in \mathcal{M}} e^{-x_m} e^{-\frac{x}{\sigma_m^2}} + 2 \sum_{j \geq m^*} e^{-x_{j,m^*}} e^{-\frac{x}{\sigma_{j,m^*}^2}}. \quad (5.12)$$

We have

$$\mathbb{E}(X^p) = \mathbb{E}(X^p \mathbb{I}_{X \geq Y + C_{m^*}}) + \mathbb{E}(X^p \mathbb{I}_{X < Y + C_{m^*}})$$
$$\leq \mathbb{E}\left[(X - Y - C_{m^*} + Y + C_{m^*})^p \mathbb{I}_{X \geq Y + C_{m^*}}\right] + \mathbb{E}\left[(Y + C_{m^*})^p \mathbb{I}_{X < Y + C_{m^*}}\right]$$
$$\leq 2^{p-1} \mathbb{E}\left[(X - Y - C_{m^*})^p \mathbb{I}_{X \geq Y + C_{m^*}}\right] + 2^{p-1} \mathbb{E}\left[(Y + C_{m^*})^p\right].$$

We have used the inequality $(a + b)^p \leq 2^{p-1}(a^p + b^p)$ which holds for all $p \geq 1$, $a, b \geq 0$.

Moreover, setting $(u)_+^p = (\max(u, 0))^p$,

$$\mathbb{E}\left[(X - Y - C_{m^*})_+^p\right] = \int_0^\infty \mathbb{P}\left((X - Y - C_{m^*})_+^p > t\right) dt$$
$$= \int_0^\infty \mathbb{P}\left((X - Y - C_{m^*})_+^p > (2\sqrt{2x})^p\right) p 2^{p-1} (\sqrt{2})^p (\sqrt{x})^{p-2} dx.$$

Hence, using (5.12), we get

$$\mathbb{E}\left((X - Y - C_{m^*})_+^p\right)$$
$$C(p) \left(\sum_{m \in \mathcal{M}} e^{-x_m} \sigma_m^p + \sum_{j \geq m^*} e^{-x_{j,m^*}} \sigma_{j,m^*}^p\right) \int_0^\infty e^{-x} (\sqrt{x})^{p-2} dx.$$

We conclude that for all $p \geq 1$, there exists some constant $C(p) > 0$ such that

$$\mathbb{E}\left[|T(\hat{s}_{\hat{m}}) - T(s)|^p\right] \leq C(p) \{\mathrm{Crit}(m^*)^p + \mathbb{E}\left(|T(\hat{s}_{m^*}) - T(s)|^p\right)\}$$
$$+ C(p) \left(\sup_{j \leq m^*} H^p(m^*, j) + \sup_{j \geq m^*} |T(s_j) - T(s)|^p\right)$$
$$+ C(p) \left(\sum_{m \in \mathcal{M}} e^{-x_m} \sigma_m^p + \sum_{j \geq m^*} e^{-x_{j,m^*}} \sigma_{j,m^*}^p + \frac{1}{n^p}\right).$$

Moreover, possibly enlarging $C(p)$,

$$\mathbb{E}\left[|T(\hat{s}_{m^*}) - T(s)|^p\right] \leq C(p) \left[|T(s_{m^*}) - T(s)|^p + \sigma_{m^*}^p\right].$$

Since

$$|T(s_j) - T(s)| \leq |T(s_j) - T(s_{m^*})| + |T(s_{m^*}) - T(s)|,$$

we obtain

$$\sup_{j \geq m^*} |T(s_j) - T(s)|^p \leq C(p) \left(|T(s_{m^*}) - T(s)|^p + \sup_{j \geq m^*} |T(s_j) - T(s_{m^*})|^p\right)$$
$$\leq C(p) \left(|T(s_{m^*}) - T(s)|^p + \mathrm{Crit}(m^*)^p\right).$$

This concludes the proof of Theorem 1.

### 5.2. Proof of Corollary 1

The proof follows from bounding the terms on the right hand side of (2.4). In the following, $C$ denotes a positive constant which may vary from line to line. We mention the dependency of these constants with respect to the parameters involved in the problem.

Let $W_j$ be the orthogonal complement of $S_j$ in $S_{j+1}$ : $S_{j+1} = W_j \bigoplus S_j$. Set $D_j(s)$ the projection of $s$ onto $W_j$.

It follows from [24, Theorem 3] p. 31 that there exists a constant $C(r)$ such that
$$\|(D_j(s))^{(r)}\|_\infty \leq C(r)2^{jr}\|D_j(s)\|_\infty .$$

We recall that
$$D_j(s) = \sum_{k\in\mathbb{Z}} \beta_{j,k}\psi_{j,k}$$
where $\psi_{j,k} = 2^{j/2}\psi(2^j x - k)$ and $\beta_{j,k} = \langle s\psi_{j,k}\rangle$.

Hence, since we have assumed that $\psi$ has a compact support,
$$\|D_j(s)\|_\infty \leq C(r)2^{j/2} \sup_{k\in\mathbb{Z}}|\beta_{j,k}|$$
and
$$\|(D_j(s))^{(r)}\|_\infty \leq C 2^{j/2+jr} \sup_{k\in\mathbb{Z}}|\beta_{j,k}|$$

Since $s \in B_{\alpha,\infty,q}([0,1])$ with $\|s\|_{\alpha,\infty,q} \leq L$,
$$\sum_{j\geq 0} 2^{qj(\alpha+1/2)} \sup_{k\in\mathbb{Z}}|\beta_{j,k}|^q \leq L^q. \tag{5.13}$$

We define $B(m) = \sup_{m\leq j\leq d_n}|s_j^{(r)}(x_0) - s_m^{(r)}(x_0)|$.

$$\begin{aligned} B(m) &\leq \sup_{j\geq m} \|s_j^{(r)} - s_m^{(r)}\|_\infty \\ &\leq \sup_{j\geq m} \sum_{l=m+1}^{j} \|(D_l(s))^{(r)}\|_\infty \\ &\leq C(r) \sum_{j>m} 2^{j(\alpha+1/2)} \sup_{k\in\mathbb{Z}}|\beta_{j,k}|2^{-j(\alpha-r)}. \end{aligned}$$

Using Cauchy-Schwarz inequality,
$$B(m) \leq C(r)\left(\sum_{j>m} 2^{qj(\alpha+1/2)}\left(\sup_{k\in\mathbb{Z}}|\beta_{j,k}|\right)^q\right)^{1/q}\left(\sum_{j>m} 2^{-jq'(\alpha-r)}\right)^{1/q'}$$
where $1/q + 1/q' = 1$. It follows from (5.13) that
$$B(m) \leq C(\alpha,r,q)L 2^{-m(\alpha-r)}.$$

Note that
$$\sigma_m^2 = \frac{\sigma^2}{n} \sum_{k \in \mathbb{Z}} \left( \varphi_{m,k}^{(r)}(x_0) \right)^2,$$

$$\forall j > m, \sigma_{j,m}^2 = \mathrm{Var}(T(\hat{s}_j) - T(\hat{s}_m)) = \frac{\sigma^2}{n} \sum_{l=m}^{j-1} \sum_{k \in \mathbb{Z}} \left( \psi_{l,k}^{(r)}(x_0) \right)^2.$$

It follows from conditions (**i**) and (**iv**) that
$$\sigma_m^2 \leq C(r) \frac{\sigma^2}{n} 2^{m(1+2r)}$$

and that
$$\sigma_{j,m}^2 \leq C(r) \frac{\sigma^2}{n} \left( 2^{j(1+2r)} - 2^{m(1+2r)} \right).$$

Thus, for all $m$ in $\mathcal{M}$,
$$\mathrm{Crit}^p(m) \leq C(p, \alpha, r, q) \left( L^p 2^{-pm(\alpha-r)} + \sigma^p \left( \frac{2^{m(1+2r)} m}{n} \right)^{p/2} \right).$$

Let
$$m_0(n) = \left[ \frac{1}{\ln(2)} \ln \left( \left( \frac{nL^2}{\sigma^2 \ln(n)} \right)^{\frac{1}{1+2\alpha}} \right) \right], \quad (5.14)$$

where $[x]$ denotes the integer part of $x$.

One can easily show that if $n/(\sigma^2 \ln(n)) \geq 2^{1+2\alpha} L^{-2}$ and $n^{2\alpha} \sigma^2 \ln(n) \geq L^2$, then
$$m_0(n) \in \mathcal{M}.$$

Hence,
$$\mathrm{Crit}(m^*) \leq \mathrm{Crit}(m_0(n)) + \frac{1}{n},$$

which implies that
$$\mathrm{Crit}^p(m^*) \leq C(p, \alpha, r, q) \left( L^{\frac{p(1+2r)}{1+2\alpha}} \left( \frac{\sigma^2 \ln n}{n} \right)^{\frac{p(\alpha-r)}{1+2\alpha}} + \frac{1}{n^p} \right).$$

We next bound the other terms on the right hand side of Theorem 1.

- For all $1 \leq m \leq d_n$,
$$|s_m^{(r)}(x_0) - s^{(r)}(x_0)|^p \leq \left( |s_{d_n}^{(r)}(x_0) - s^{(r)}(x_0)| + |s_m^{(r)}(x_0) - s_{d_n}^{(r)}(x_0)| \right)^p$$

$$\leq C(\alpha, r, q) \left( L^p 2^{-pd_n(\alpha-r)} + \sup_{j \geq m, j \in \mathcal{M}} |s_m^{(r)}(x_0)) - s_j^{(r)}(x_0)|^p \right).$$

This implies that
$$|s_{m^*}^{(r)}(x_0) - s^{(r)}(x_0)|^p \leq C(\alpha, r, q) \left( L^p 2^{-pd_n(\alpha-r)} + \mathrm{Crit}^p(m^*) \right).$$

- For all $1 \leq m \leq d_n$, $\sigma_m \leq pen(m) \leq \text{Crit}(m)$.
- For all $1 \leq m \leq d_n$,

$$\sup_{j \leq m} H(m,j) = \sup_{j \leq m} \sqrt{x_{m,j}} \sigma_{m,j} \leq \sqrt{x_m} \sigma_m \leq \text{Crit}(m) .$$

This implies that

$$\text{Crit}^p(m^*) + |s_{m^*}^{(r)}(x_0) - s^{(r)}(x_0)|^p + \sigma_{m^*}^p + \sup_{j \leq m^*} H^p(j,m^*)$$

$$\leq C(p,\alpha,r,q) \left( L^{\frac{p(1+2r)}{1+2\alpha}} \left( \frac{\sigma^2 \ln n}{n} \right)^{\frac{p(\alpha-r)}{1+2\alpha}} + \frac{1}{n^p} + L^p 2^{-pd_n(\alpha-r)} \right) .$$

Since $2^{-d_n} < 2/n$ and $n^{2\alpha} \sigma^2 \ln n \geq L^2$,

$$L^{\frac{p(1+2r)}{1+2\alpha}} \left( \frac{\sigma^2 \ln n}{n} \right)^{\frac{p(\alpha-r)}{1+2\alpha}} \geq L^p 2^{-pd_n(\alpha-r)} 2^{p(\alpha-r)}.$$

On the other hand we have :

$$\sum_{m \in \mathcal{M}} e^{-x_m} \sigma_m^p + \sum_{j \geq m^*} e^{-x_{j,m^*}} \sigma_{j,m^*}^p \leq C(r,p) \frac{\sigma^p}{n^{p/2}} d_n \leq C(r,p) \sigma^p \frac{\ln(n)}{n^{p/2}} ,$$

which yields the desired bound.

### 5.3. Proof of Corollary 2

As in Corollary 1, the proof follows by bounding the terms on the right hand side of (2.4). Let us first control $\sigma_m^2$ for all $m \in \mathcal{M}$. For $m \leq m_n$,

$$\sigma_m^2 = \frac{\sigma^2}{n} \sum_{k \in \mathbb{Z}} T^2(\phi_{m,k}) = \frac{\sigma^2}{nH_n^2} \sum_{k \in \mathbb{Z}} 2^m \int_{I_{H_n}} \varphi(2^m x - k) \int_{I_{H_n}} \varphi(2^m x - k).$$

Since $\varphi$ has a compact support,

$$\sum_{k \in \mathbb{Z}} \left| \int_{I_{H_n}} \varphi(2^m x - k) \right| \leq C \|\varphi\|_\infty H_n,$$

and

$$\sup_{k \in \mathbb{Z}} \left| \int_{I_{H_n}} \varphi(2^m x - k) \right| \leq C \|\varphi\|_\infty (2^{-m} \wedge H_n).$$

Hence, for all $m \leq m_n$,

$$\sigma_m^2 \leq C \frac{\sigma^2}{n} 2^m.$$

Moreover,

$$\sigma^2_{m_n+1} = \frac{\sigma^2}{n} T^2(\sqrt{H_n}\, \mathbb{1}_{I_{H_n}}) \leq \sigma^2/(nH_n) \leq \sigma^2 2^{m_n+1}/n.$$

It is easy to see that $\sigma^2_{j,m} \leq 2(\sigma^2_m + \sigma^2_j)$ and that $\sigma_{m,m} = 0$. Hence

$$\forall j \neq m \in \mathcal{M}, \sigma^2_{j,m} \leq C\frac{\sigma^2}{n} 2^{m \vee j}.$$

This implies that

$$\sum_{m \in \mathcal{M}} e^{-x_m} \sigma^p_m \leq C(p) \frac{\sigma^p}{n^{p/2}}$$

$$\forall m \in \mathcal{M}, \sum_{j \geq m} e^{-x_{j,m}} \sigma^p_{j,m} \leq C(p) \frac{\sigma^p}{n^{p/2}}.$$

Since $T(s_{m_n+1}) = T(s)$, and since for all $f, g \in \mathbb{L}^2([0,1])$, $|T(f) - T(g)| \leq \|f - g\|_\infty$, we obtain, by the same computations as in the proof of Corollary 1, that for all $m \leq m_n$

$$\sup_{j > m, j \in \mathcal{M}} |T(s_m) - T(s_j)| \leq \sup_{m < j \leq m_n} \|s_m - s_j\|_\infty + \|s_m - s\|_\infty \leq C(\alpha, q) L 2^{-m\alpha}.$$

Let $m_0(n)$ be defined by (5.14). Since we assumed that $nL^2/\sigma^2 \ln(n) \geq 1$, $m_0(n) \in \mathbb{N}$. If $\left(nL^2/\sigma^2 \ln(n)\right)^{1/(1+2\alpha)} H_n \leq 1$, $2^{m_0(n)} \leq 1/H_n$ and we get

$$\begin{aligned}
\mathrm{Crit}^p(m^*) &\leq C(p)\left(\mathrm{Crit}^p(m_0(n)) + n^{-p}\right) \\
&\leq C(p, \alpha, q) \left( L^{\frac{p}{1+2\alpha}} \left(\frac{\sigma^2 \ln n}{n}\right)^{\frac{p\alpha}{2\alpha+1}} + n^{-p}\right).
\end{aligned}$$

Moreover,

$$\begin{aligned}
\mathrm{Crit}^p(m^*) &\leq C(p)\left(\mathrm{Crit}^p(m_n+1) + n^{-p}\right) \\
&\leq C(p) \left( \frac{\sigma^p}{(nH_n)^{p/2}} \left(\ln\left(\frac{1}{H_n}\right)\right)^{p/2} + n^{-p}\right).
\end{aligned}$$

Also, since $T(s_{m_n+1}) = T(s)$,

$$|T(s_{m^*}) - T(s)| \leq \sup_{j \geq m^*, j \in \mathcal{M}} |T(s_{m^*}) - T(s_j)| \leq \mathrm{Crit}(m^*).$$

The other terms appearing in the upper bound for $E(|T(\hat{s}_{\hat{m}}) - T(s)|^p)$ given in Theorem 1 can be controled as in the proof of Corollary 1.

## 5.4. Proof of Lemma 1

It follows from the results of Donoho and Liu [13] that

$$\inf_{T_n} \sup_{s \in \mathcal{H}_\alpha(L)} \mathbb{E}\left((T_n - T(s))^2\right) \geq C(\sigma)\omega_2^2\left(\frac{1}{\sqrt{n}}, T, \mathcal{H}_\alpha(L)\right),$$

where $\omega_2(\epsilon, T, \mathcal{F})$ denotes the modulus of continuity of the linear functional $T$ over the set $\mathcal{F}$ with respect to the $\mathbb{L}^2$ norm, namely

$$\omega_2(\epsilon, T, \mathcal{F}) = \sup\left\{|T(f_1) - T(f_0)|, f_0, f_1 \in \mathcal{F}, \|f_1 - f_0\|_2 \leq \epsilon\right\}.$$

In order to prove (3.6), we consider the functions $f_0 = 0$ and $f_1$ defined on $[0,1]$ by:

$$\begin{aligned} f_1(x) &= \rho_n x^\alpha \text{ for } x \in [0, H_n], \\ f_1(x) &= \rho_n(2H_n - x)^\alpha \text{ for } x \in (H_n, 2H_n], \\ f_1(x) &= 0 \text{ for } x > 2H_n \end{aligned}$$

where

$$\rho_n = ((1 + 2\alpha)/(2nH_n^{1+2\alpha}))^{1/2}.$$

The condition $H_n \geq \left((1+2\alpha)/(2nL^2)\right)^{1/(1+2\alpha)}$ ensures that $\rho_n \leq L$, hence $f_1 \in \mathcal{H}_\alpha(L)$. One can easily verify that $\|f_1 - f_0\|_2 = 1/\sqrt{n}$ and that $T(f_1) = C(\alpha)/\sqrt{nH_n}$.

In order to prove (3.7), we consider the functions $f_0 = 0$ and $f_1$ defined on $[0,1]$ by:

$$\begin{aligned} f_1(x) &= L(\gamma_n - x)^\alpha \text{ for } x \in [0, \gamma_n], \\ f_1(x) &= 0 \text{ for } x > \gamma_n \end{aligned}$$

where $\gamma_n = ((1+2\alpha)/(nL^2))^{1/(1+2\alpha)}$. Since we have assumed that $nL^2 \geq 1+2\alpha$, $\gamma_n \leq 1$. We have, $\|f_1 - f_0\|_2 = 1/\sqrt{n}$ and

$$T(f_1) = \frac{L\gamma_n^{1+\alpha}}{H_n(1+\alpha)}\left[1 - \left(1 - \frac{H_n}{\gamma_n}\right)^{1+\alpha}\right].$$

$$T(f_1) = L\gamma_n^\alpha C_n^\alpha$$

where $C_n \in [1 - H_n/\gamma_n, 1]$. The last equality is obtained by Taylor-Lagrange's formula. For $H_n \leq \frac{1}{2}(1+2\alpha)^{1/(1+2\alpha)}\left(nL^2\right)^{-1/(1+2\alpha)}$, $H_n/\gamma_n \leq 1/2$, which leads to $T(f_1) \geq C(\alpha)L\gamma_n^\alpha$, hence (3.7) holds.

## 5.5. Proof of Corollary 3

For all $m \in \mathcal{M}$, we set $\Lambda_m = \left\{0, 1, \ldots, 2^{m-1}\right\}^d$, and for all $\lambda = (k_1, \ldots, k_d) \in \Lambda_m$, let $I_\lambda = [\frac{k_1}{2^m}, \frac{k_1+1}{2^m}[\times \cdots \times [\frac{k_d}{2^m}, \frac{k_d+1}{2^m}[$ and

$$\phi_\lambda = 2^{md/2}\, \mathbb{1}_{I_\lambda}.$$

Then $(\phi_\lambda, \lambda \in \Lambda_m)$ is an orthonormal basis of $S_m$. It is easy to verify that for all $m \in \mathcal{M}$,
$$\sigma_m^2 = \frac{\sigma^2}{n} 2^{md}$$
and that for all $(j,m) \in \mathcal{M}^2$,
$$\sigma_{j,m}^2 \leq 2\frac{\sigma^2}{n}(2^{md} + 2^{jd}).$$
It follows from the definition of $x_m$ and $x_{j,m}$ that
$$\sum_{m \in \mathcal{M}} e^{-x_m} \sigma_m^p \leq C \frac{\sigma^p}{n^{p/2}} d_n,$$
$$\sum_{j \geq m^*} e^{-x_{j,m^*}} \sigma_{j,m^*}^p \leq \frac{\sigma^p}{n^{p/2}} d_n.$$
Moreover,
$$\sigma_{m^*} + \sup_{j \leq m^*} \sqrt{x_{m^*,j}} \sigma_{m^*,j} \leq C\mathrm{pen}(m^*)$$
and
$$|s_{m^*}(x_0) - s(x_0)|^p \leq C(p) \sup_{j \geq m^*, j \in \mathcal{M}} |s_j(x_0) - s_{m^*}(x_0)|^p$$
$$+ C(p)|s_{d_n}(x_0) - s(x_0)|^p.$$
Hence, it follows from Theorem 1 that, possibly enlarging $C(\epsilon)$,
$$\mathbb{E}\left(|\hat{s}_{\hat{m}}(x_0) - s(x_0)|^p\right) \leq C(p)\left(\mathrm{Crit}^p(m^*) + \frac{\sigma^p d_n}{n^{p/2}} + |s_{d_n}(x_0) - s(x_0)|^p\right).$$
For all $m \in \mathcal{M}$,
$$s_m = \sum_{\lambda \in \Lambda_m} \frac{1}{2^{md}}\left(\int_{I_\lambda} s(x)dx\right) \mathbb{1}_{I_\lambda}.$$
This implies that,
$$s_m(x_0) - s(x_0) = \frac{1}{2^{md}}\int_{I_\lambda(x_0)} (s(x) - s(x_0))dx,$$
where $I_\lambda(x_0)$ is the set $I_\lambda$ that contains $x_0$. Hence, if $s \in \mathcal{H}_\alpha(L)$,
$$|s_m(x_0) - s(x_0)| \leq L 2^{-m\alpha}.$$
Let
$$m_1(n) = \left[\frac{1}{\ln(2)}\ln\left(\left(\frac{nL^2}{\ln(n)}\right)^{\frac{1}{d+2\alpha}}\right)\right].$$

Since $m_1(n) \in \mathcal{M}$, as soon as $n/\ln(n) \geq 2^{d+2\alpha} L^{-2}$ and $L^2 \leq n^{\frac{2\alpha}{d}} \ln(n)$,

$$Crit^p(m^*) \leq Crit^p(m_1(n)) \leq C(p, \alpha, d, \sigma) L^{\frac{pd}{d+2\alpha}} \left(\frac{\ln(n)}{n}\right)^{\frac{p\alpha}{d+2\alpha}}.$$

Hence, for $n$ large enough,

$$\mathbb{E}\left(|\hat{s}_{\hat{m}}(x_0) - s(x_0)|^p\right) \leq C(p, \alpha, d, \sigma) L^{\frac{pd}{d+2\alpha}} \left(\frac{\ln(n)}{n}\right)^{\frac{p\alpha}{d+2\alpha}}.$$

This concludes the proof of Corollary 3.

**Acknowledgment**

We would like to thank the Associate Editor, the referees and A. Tsybakov for their helpful comments.